\documentclass[a4paper]{article}

\usepackage{arydshln}
\usepackage{graphicx}
\usepackage{amsmath}
\usepackage{amssymb, verbatim}
\usepackage{mathrsfs}
\usepackage{psfig}
\usepackage{dsfont}
\usepackage{url}
\usepackage[usenames,dvipsnames]{color}
\usepackage[compress]{cite}
\usepackage{mdwlist}
\usepackage{paralist}
\usepackage[hyperindex,breaklinks]{hyperref}

\newtheorem{lemma}{\sc Lemma}[section]
\newtheorem{theorem}[lemma]{\sc Theorem}

\newtheorem{definition}{\sc Definition}[section]

\newcommand{\jpfig}[4]{\begin{figure}[t] \centering \includegraphics[width=#1\linewidth]{#2} \caption{\label{#3}#4} \end{figure}}
\renewcommand{\matrix}[2]{\left[\begin{array}{#1} #2 \end{array}\right] }
\newcommand{\diag}{\text{diag}}
\newcommand{\p}{\mathcal{P}}
\newcommand{\A}{\mathcal{A}}
\newcommand{\B}{\mathcal{B}}
\newcommand{\D}{\mathcal{D}}
\newcommand{\R}{\mathcal{R}}
\newcommand{\K}{\mathcal{K}}
\newcommand{\comm}{\mathcal{C}}

\def\IEEEQED{~\rule[-1pt]{5pt}{5pt}\par\medskip}
\newenvironment{IEEEproof}{{\it Proof:\ }}{ \hfill \IEEEQED}

\begin{document}

\title{\bf Optimal Disturbance Accommodation with Limited Model Information\thanks{ The work of F.~Farokhi and K.~H.~Johansson was supported by grants from the Swedish Research Council and the Knut and Alice Wallenberg Foundation. The work of C.~Langbort was supported, in part, by the 2010 AFOSR MURI ``Multi-Layer and Multi-Resolution Networks of Interacting Agents in Adversarial Environments''.}}

\author{Farhad~Farokhi\thanks{F. Farokhi and K. H. Johansson are with ACCESS Linnaeus Center, School of Electrical Engineering, KTH-Royal Institute of Technology, SE-100 44 Stockholm, Sweden. E-mails: \{farokhi,kallej\}@ee.kth.se },~C\'{e}dric~Langbort\thanks{C. Langbort is with the Department of Aerospace Engineering, University of Illinois at Urbana-Champaign, Illinois, USA. E-mail: langbort@illinois.edu},~and~Karl~H.~Johansson$^\dag$}

\maketitle

\begin{abstract}
The design of optimal dynamic disturbance accommodation controller with limited model information is considered. We adapt the family of limited model information control design strategies, defined earlier by the authors, to handle dynamic controllers. This family of limited model information design strategies construct subcontrollers distributively by accessing only local plant model information. The closed-loop performance of the dynamic controllers that they can produce are studied using a performance metric called the competitive ratio which is the worst case ratio of the cost a control design strategy to the cost of the optimal control design with full model information.
\end{abstract}

\section{Introduction}
Recent advances in networked control engineering have opened new doors toward controlling large-scale systems. These large-scale systems are naturally composed of many smaller unit that are coupled to each other~\cite{Giulietti2000,Swaroop1999,Dunbar2007,Negenborn2010}. For these large-scale interconnected systems, we can either design a centralized or a decentralized controller. Contrary to a centralized controller, each subcontroller in a decentralized controller only observes a local subset of the state-measurements~(e.g.,~\cite{Levine1971,Rotkowitz2006,Voulgaris200351}). When designing these controllers, generally, it is assumed that the global model of the system is available to each subcontroller's designer. However, there are several reasons why such plant model information would not be globally known. One reason could be that the subsystems consider their model information private, and therefore, they are reluctant to share information with other subsystems. This case can be well illustrated by supply chains or power networks where the economic incentives of competing companies might limit the exchange of model information between the companies. It might also be the case that the full model is not available at the moment, or the designer would like to not modify a particular subcontroller, if the model of a subsystem changes. For instance, in the case of cooperative driving, each vehicle controller simply cannot be designed based on model information of all possible vehicles that it may interact with in future. Therefore, we are interested in finding control design strategies which construct subcontrollers distributively for plants made of interconnected subsystems without the global model of the system. The interconnection structure and the common closed-loop cost to be minimized are assumed to be public knowledge. We identify these control design methods by ``limited model information'' control design strategies~\cite{Farokhi_ACC_2011,Farokhi_Technical_2011}.

Multi-variable servomechanism and disturbance accommodation control design is one of the oldest problems in control engineering~\cite{Johnson1968}. We adapt the procedure introduced in~\cite{Johnson1968,Anderson1971} to design optimal disturbance accommodation controllers for discrete-time linear time-invariant plants under a separable quadratic performance measure. The choice of the cost function is motivated first by the optimal disturbance accommodation literature~\cite{Johnson1968,Anderson1971}, and second by our interest in dynamically-coupled but cost-decoupled plants and their applications in supply chains and shared infrastructures~\cite{Dunbar2007,Negenborn2010}. Then, we focus on the disturbance accommodation design problem under limited model information. We investigate the achievable closed-loop performance of the dynamic controllers that the limited model information control design strategies can produce using the competitive ratio, that is, the worst case ratio of the cost a control design strategy to the cost of the optimal control design with full model information. We find a minimizer of the competitive ratio over the set of limited model information control design strategies. Since this minimizer may not be unique we prove that it is undominated, that is, there is no other control design method that acts better while exhibiting the same worst-case ratio.

This paper is organized as follows. We mathematically formulate the problem in Section~\ref{sec_2}. In Section~\ref{sec_2_1}, we introduce two useful control design strategies and study their properties. We characterize the best limited model information control design method as a function of the subsystems interconnection pattern in Section~\ref{sec_3}. In Section~\ref{sec:Gc}, we study the trade-off between the amount of the information available to each subsystem and the quality of the controllers that they can produce. Finally, we end with conclusions in Section~\ref{sec_4}.

\subsection{Notation}
The set of real numbers and complex numbers are denoted by $\mathbb{R}$ and $\mathbb{C}$, respectively. All other sets are denoted by calligraphic letters such as $\p$ and $\A$. The notation $\R$ denotes the set of proper real rational functions.

Matrices are denoted by capital roman letters such as $A$. $A_j$ will denote the $j^{\textrm{th}}$ row of $A$. $A_{ij}$ denotes a sub-matrix of matrix $A$, the dimension and the position of which will be defined in the text. The entry in the $i^{\textrm{th}}$ row and the $j^{\textrm{th}}$ column of the matrix $A$ is $a_{ij}$.

Let $\mathcal{S}_{++}^n$ ($\mathcal{S}_{+}^n$) be the set of symmetric positive definite (positive semidefinite) matrices in $\mathbb{R}^{n\times n}$. $A > (\geq) 0$ means that the symmetric matrix $A\in \mathbb{R}^{n\times n}$ is positive definite (positive semidefinite) and $A > (\geq) B$ means $A-B > (\geq) 0$.

$\underline{\sigma}(Y)$ and $\overline{\sigma}(Y)$ denote the smallest and the largest singular values of the matrix $Y$, respectively. Vector $e_i$ denotes the column-vector with all entries zero except the $i^{\textrm{th}}$ entry, which is equal to one.

All graphs considered in this paper are directed, possibly with self-loops, with vertex set $\{1,\dots,n\}$ for some positive integer $n$. We say $i$ is a sink in $G=(\{1,\dots,n\},E)$, if there does not exist $j \neq i $ such that $(i,j) \in E$. The adjacency matrix $S\in\{0,1\}^{n\times n}$ of graph $G$ is a matrix whose entries are defined as $s_{ij}=1$ if $(j,i) \in E$ and $s_{ij}=0$ otherwise. Since the set of vertices is fixed here, a subgraph of $G$ is a graph whose edge set is a subset of the edge set of $G$ and a supergraph of $G$ is a graph of which $G$ is a subgraph. We use the notation $G'\supseteq G$ to indicate that $G'$ is a supergraph of $G$.

\section{Mathematical Formulation} \label{sec_2}
\subsection{Plant Model}
Consider the discrete-time linear time-invariant dynamical system described in state-space representation by
\begin{equation} \label{eqn_1}
x(k+1)=Ax(k)+B(u(k)+w(k)) \; ; \; x(0)=x_0,
\end{equation}
where $x(k)\in \mathbb{R}^{n}$ is the state vector, $u(k)\in \mathbb{R}^{n}$ is the control input, and $w(k)\in \mathbb{R}^{n}$ is the disturbance vector. In addition, assume that $w(k)$ is a dynamic disturbance modeled as
\begin{equation} \label{eqn_2}
w(k+1)=Dw(k) \; ; \; w(0)=w_0.
\end{equation}
Let a plant graph $G_\p$ with adjacency matrix $S_\p$ be given. We define the following set of matrices
\begin{equation*}
\begin{split}
\A(S_{\p})=\{ \bar{A} \in \mathbb{R}^{n \times n} \; | \; \bar{a}_{ij} = 0  & \mbox{ for all }  1\leq i,j \leq n \mbox { such that } (s_{\p})_{ij}=0 \}.
\end{split}
\end{equation*}
Also, let us define
\begin{equation*}
\begin{split}
\B(\epsilon_b) =\{ \bar{B} \in  \mathbb{R}^{n \times n} \; | \; \underline{\sigma}(\bar{B}) &\geq \epsilon_b, \bar{b}_{ij} = 0 \mbox{ for all } 1 \leq i \neq j \leq n \},
\end{split}
\end{equation*}
for a given scalar $\epsilon_b >0$ and
\begin{equation*}
\D=\{ \bar{D} \in \mathbb{R}^{n \times n} \; | \; \bar{d}_{ij} = 0 \mbox{ for all } 1 \leq i \neq j \leq n \}.
\end{equation*}
We can introduce the set of plants of interest $\p$ as the space of all discrete-time linear time-invariant systems of the form~(\ref{eqn_1})~and~(\ref{eqn_2}) with $A \in \A (S_\p)$, $B \in \B (\epsilon_b)$, $D \in \D$, $x_0 \in \mathbb{R}^n$, and $w_0 \in \mathbb{R}^n$. Since $\p$ is isomorph to $\A(S_{\p}) \times \B(\epsilon_b) \times \D \times \mathbb{R}^n \times \mathbb{R}^n$, we identify a plant $P \in \p$ with its corresponding tuple $(A,B,D,x_0,w_0)$ with a slight abuse of notation.

We can think of $x_i\in\mathbb{R}$, $u_i\in\mathbb{R}$, and $w_i\in\mathbb{R}$ as the state, input, and disturbance of scalar subsystem $i$ with its dynamic given as
$$
x_i(k+1) = \sum_{j=1}^n a_{ij} x_j(k) + b_{ii} (u_i(k)+w_i(k)).
$$
We call $G_\p$ the plant graph since it illustrates the interconnection structure between different subsystems, that is, subsystem $j$ can affect subsystem $i$ only if $(j,i) \in E_{\p}$. In this paper, we assume that overall system is fully-actuated, that is, any $B\in\B(\epsilon_b)$ is a square invertible matrix. This assumption is motivated by the fact that we want all the subsystems to be directly controllable.

\subsection{Controller} \label{subsec:controller}
The control laws of interest in this paper are discrete-time linear time-invariant dynamic state-feedback control laws of the form
$$
x_K(k+1)=A_K x_K(k) + B_K x(k) \; ; \; x_K(0)=0,
$$
$$
u(k)=C_K x_K(k)+D_K x(k).
$$
Each controller can also be represented by its transfer function
$$
K\triangleq\left[\begin{array}{c|c} A_K & B_K \\ \hline C_K & D_K \end{array}\right]=C_K(zI-A_K)^{-1}B_K+D_K,
$$
where $z$ is the symbol for one time-step forward shift operator. Let a control graph $G_{\K}$ with adjacency matrix $S_{\K}$ be given. Each controller $K$ must belong to
\begin{equation*}
\begin{split}
\K(S_{\K})=\{ \bar{K} \in  \; \R^{n \times n} \;|\; \bar{k}_{ij} = 0 &\mbox{ for all }  1\leq i,j \leq n \mbox { such that } (s_{\K})_{ij}=0 \}.
\end{split}
\end{equation*}
When adjacency matrix $S_{\K}$ is not relevant or can be deduced from context, we refer to the set of controllers as $\K$. Since it makes sense for each subsystem's controller to have access to at least its own state-measurements, we make the standing assumption that in each control graph $G_\K$, all the self-loops are present.

Finding the optimal structured controller is difficult (numerically intractable) for general $G_\K$ and $G_\p$ even when the global model is known. Therefore, in this paper, as a starting point, we only concentrate on the cases where the control graph $G_{\K}$ is a supergraph of the plant graph $G_\p$.

\subsection{Control Design Methods}
A control design method $\Gamma$ is a mapping from the set of plants $\p$ to the set of controllers $\K$. We can write the control design method $\Gamma$ as
\begin{equation*}
\Gamma=\matrix{ccc}{ \gamma_{11} & \cdots & \gamma_{1n} \\ \vdots & \ddots & \vdots \\ \gamma_{n1} & \cdots & \gamma_{nn} }
\end{equation*}
where each entry $\gamma_{ij}$ represents a map $\A(S_{\p}) \times \B(\epsilon_b) \times \D  \rightarrow \R$. Let a design graph $G_{\mathcal{C}}$ with adjacency matrix $S_{\mathcal{C}}$ be given. The control design strategy $\Gamma$ has structure $G_{\mathcal{C}}$ if, for all $i$, the map $\Gamma_i=[\gamma_{i1}\;\cdots\;\gamma_{in}]$ is only a function of $\left\{ [a_{j1}\;\cdots\;a_{jn}] , b_{jj}, d_{jj} \; | \; (s_{\mathcal{C}})_{ij} \neq 0 \right\}$. Consequently, for each $i$, subcontroller $i$ is constructed with model information of only those subsystems $j$ that $(j,i) \in E_{\mathcal{C}}$. We are only interested in those control design strategies that are neither a function of the initial state $x_0$ nor of the initial disturbance~$w_0$. The set of all control design strategies with the design graph $G_{\mathcal{C}}$ is denoted by $\mathcal{C}$. Since it makes sense for the designer of each subsystem's controller to have access to at least its own model parameters, we make the standing assumption that in each design graph $G_\comm$, all the self-loops are present.

For simplicity of notation, let us assume that any control design strategy $\Gamma\in\comm$ has a state-space realization of the form
$$
\Gamma(A,B,D)=\left[\begin{array}{c|c} A_\Gamma(A,B,D) & B_\Gamma(A,B,D) \\ \hline C_\Gamma(A,B,D) & D_\Gamma(A,B,D) \end{array}\right],
$$
where matrices $A_\Gamma(A,B,D)$, $B_\Gamma(A,B,D)$, $C_\Gamma(A,B,D)$, and $D_\Gamma(A,B,D)$ are of appropriate dimension for each plant $P=(A,B,D,x_0,w_0)\in \p$. The matrices $A_\Gamma(A,B,D)$ and $C_\Gamma(A,B,D)$ are block diagonal matrices since different subcontrollers should not share state variables. This realization is not necessarily a minimal realization.

\subsection{Performance Metrics}
We need to introduce performance metrics to compare the control design methods. These performance metrics are adapted from earlier definitions in~\cite{Langbort2010,Farokhi_ACC_2011}. Let us start with introducing the closed-loop performance criterion.

To each plant $P=(A,B,D,x_0,w_0) \in \p$ and controller $K\in \K$, we associate the performance criterion
\begin{equation*}
J_P (K)\hspace{-.03in}=\hspace{-.03in}\sum_{k=0}^\infty [ x(k)^TQx(k) + (u(k)+w(k))^TR(u(k)+w(k))]
\end{equation*}
where $Q \in \mathcal{S}_{++}^n$ and $R \in \mathcal{S}_{++}^n$ are diagonal matrices. We make the standing assumption that $Q=R=I$. This is without loss of generality because of the change of variables $(\bar{x},\bar{u},\bar{w})= (Q^{1/2} x,R^{1/2}u,R^{1/2}w)$ that transforms the state-space representation into
\begin{equation*}
\begin{split}
\bar{x}(k+1)&\hspace{-.03in}=Q^{1/2}AQ^{-1/2}\bar{x}(k)
\hspace{-.03in}+\hspace{-.02in}Q^{1/2}BR^{-1/2}(\bar{u}(k)\hspace{-.03in}+\bar{w}(k))
\\&\hspace{-.03in}=\bar{A}\bar{x}(k)+\bar{B}(\bar{u}(k)+\bar{w}(k)),
\end{split}
\end{equation*}
and the performance criterion into
\begin{equation}
\label{cost_easy}
J_P (K)\hspace{-.03in}=\hspace{-.03in}\sum_{k=0}^\infty [\bar{x}(k)^T\bar{x}(k)+(\bar{u}(k)+\bar{w}(k))^T(\bar{u}(k)+\bar{w}(k))].
\end{equation}
This change of variable would not affect the plant, control, or design graph since both $Q$ and $R$ are diagonal matrices.

\begin{definition}\emph{(Competitive Ratio)} Let a plant graph $G_{\p}$ and a constant $\epsilon_b > 0$ be given. Assume that, for every plant $P \in \p$, there exists an optimal controller $K^*(P) \in \K$ such that
\begin{equation*}
J_P (K^*(P))\leq J_P (K), \; \forall K \in \K.
\end{equation*}
The competitive ratio of a control design method $\Gamma$ is defined as
\begin{equation*}
r_{\p} (\Gamma)= \sup_{P=(A,B,D,x_0,w_0) \in \p} \frac{J_P (\Gamma(A,B,D))}{J_P (K^*(P))},
\end{equation*}
with the convention that ``$\frac{0}{0}$'' equals one.
\label{def_comp_rat}
\end{definition}

\begin{definition}\emph{(Domination)} A control design method $\Gamma$ is said to dominate another control design method $\Gamma'$ if for all plants $P=(A,B,D,x_0,w_0)\in \p$
\begin{equation}
J_P(\Gamma(A,B,D))\leq J_P(\Gamma'(A,B,D)),
\label{comp}
\end{equation}
with strict inequality holding for at least one plant in $\p$. When $\Gamma' \in \mathcal{C}$ and no control design method $\Gamma \in \mathcal{C}$ exists that dominates it, we say that $\Gamma'$ is undominated in $\mathcal{C}$.
\end{definition}

\subsection{Problem Formulation}
For a given plant graph $G_\p$, control graph $G_\K$, and design graph $G_\comm$, we want to solve the problem
\begin{equation} \label{eqn:0}
\arg \min_{\Gamma \in \comm} r_{\p} (\Gamma).
\end{equation}
Because the solution to this problem might not be unique, we also want to determine which ones of these minimizers are undominated.

\section{Preliminary Results} \label{sec_2_1}
In order to give the main results of the paper, we need to introduce two control design strategies and study their properties.

\subsection{Optimal Centralized Control Design Strategy}
In this subsection, we find the optimal centralized control design strategy $K^*_C(P)$ for all plants $P\in \p$; i.e., the optimal control design strategy when the control graph $G_\K$ is a complete graph. Note that we use the notation $K^*_C(P)$ to denote the centralized optimal control design strategy as the notation $K^*(P)$ is reserved for the optimal control design strategy for a given control graph $G_\K$. We adapt the procedure given in~\cite{Anderson1971,Johnson1968} for constant input-disturbance rejection in  continuous-time systems to our framework.

First, let us define the auxiliary variables $\xi(k)=u(k)+w(k)$ and $\bar{u}(k)=u(k+1)-Du(k)$. It is evident that
\begin{equation} \label{eqn:xi}
\begin{split}
\xi(k+1)=D\xi(k)+\bar{u}(k).
\end{split}
\end{equation}
Augmenting~(\ref{eqn:xi}) with the system state-space representation in~(\ref{eqn_1}) results in
\begin{equation} \label{eqn:new_plant}
\matrix{c}{x(k+1)\\ \xi(k+1)}=\matrix{cc}{A & B \\0 & D}\matrix{c}{x(k)\\ \xi(k)}+\matrix{c}{0\\ I}\bar{u}(k).
\end{equation}
In addition, we can write the performance measure in~(\ref{cost_easy}) as
\begin{equation} \label{eqn:new_cost}
J_P (K)=\sum_{k=0}^\infty \matrix{c}{x(k)\\ \xi(k)}^T\matrix{c}{x(k)\\ \xi(k)}.
\end{equation}
To guarantee existence and uniqueness of the optimal controller $K^*_C(P)$ for any given plant $P\in\p$, we need the following lemma to hold~\cite{Molinari1975}.

\begin{lemma} \label{prop:0} The pair $(\tilde{A},\tilde{B})$ with
\begin{equation} \label{ABtilde}
\tilde{A}=\matrix{cc}{A & B \\0 & D},\hspace{.3in} \tilde{B}=\matrix{c}{0\\ I},
\end{equation}
is controllable for any given $P=(A,B,D,x_0,w_0)\in\p$. \end{lemma}

\begin{IEEEproof} The pair $(\tilde{A},\tilde{B})$ is controllable if and only if
$$
\matrix{c;{2pt/2pt}c}{\tilde{A}-\lambda I & \tilde{B}}=\matrix{cc;{2pt/2pt}c}{A-\lambda I & B & 0 \\ 0 & D-\lambda I & I }
$$
is full-rank for all $\lambda\in \mathbb{C}$. This condition is always satisfied since all the matrices $B\in\B(\epsilon_b)$ are full-rank matrices.
\end{IEEEproof}

Now, the problem of minimizing the cost function in~(\ref{eqn:new_cost}) subject to plant dynamics in~(\ref{eqn:new_plant}) becomes a state-feedback linear quadratic optimal control design with a unique solution of the form
$$
\bar{u}(k)=G_1 x(k) + G_2 \xi(k)
$$
where $G_1\in \mathbb{R}^{n\times n}$ and $G_2 \in \mathbb{R}^{n\times n}$. Therefore, we have
\begin{equation} \label{eqn:uk+11}
\begin{split}
u(k+1)&=Du(k)+\bar{u}(k)\\&=Du(k)+G_1 x(k)+G_2 \xi(k).
\end{split}
\end{equation}
Using $\xi(k)=B^{-1}(x(k+1)-Ax(k))$ in~(\ref{eqn:uk+11}), we get
\begin{equation} \label{eqn:uk+1}
\begin{split}
u(k+1)=Du(k)&+G_1 x(k) +G_2 B^{-1}(x(k+1)-Ax(k)).
\end{split}
\end{equation}
Putting a control signal of the form $u(k)=x_K(k)+D_K x(k)$ in~(\ref{eqn:uk+1}) results in
\begin{equation*}
\begin{split}
x_K(k+1)=Dx_K(k)+&(DD_K+G_1-G_2 B^{-1}A)x(k)+(G_2 B^{-1}-D_K)x(k+1).
\end{split}
\end{equation*}
Now, because of the form of the control laws of interest introduced earlier in Subsection~\ref{subsec:controller}, we have to enforce $G_2 B^{-1}-D_K=0$. Therefore, the optimal controller $K^*_C(P)$ becomes
$$
x_K(k+1)=Dx_K(k)+[G_1 +DG_2 B^{-1}-G_2 B^{-1}A]x(k),
$$
$$
u(k)=x_K(k)+G_2 B^{-1} x(k),
$$
with the initial condition $x_K(0)=0$ again because of the form of the control laws of interest.

\begin{lemma} \label{lem:1} Let the control graph $G_{\K}$ be a complete graph.  Then, the cost of the optimal control design strategy $K^*_C$ for each plant $P\in \p$ is lower-bounded as
\begin{equation*}
J_P(K^*_C(P))\geq \matrix{c}{x_0 \\ Bw_0}^T \matrix{cc}{V_{11} & V_{12} \\ V_{12}^T & V_{22} } \matrix{c}{x_0 \\ Bw_0},
\end{equation*}
where
\begin{eqnarray}
V_{11}&=&W+D^2B^{-2}+D W D, \label{eqn:V1} \\ V_{12}&=&-D(W+B^{-2}), \label{eqn:V2} \\ V_{22}&=&W+B^{-2}, \label{eqn:V3}
\end{eqnarray}
with the matrix $W$ defined as
\begin{equation} \label{eqn:W}
W=A^T(I+B^2)^{-1}A+I.
\end{equation}
\end{lemma}

\begin{IEEEproof}
To make the proof easier, let us define
\begin{equation*}
\bar{J}_P (K,\rho)=\sum_{k=0}^\infty \left( \matrix{c}{x(k)\\ \xi(k)}^T\matrix{c}{x(k)\\ \xi(k)} + \rho \bar{u}(k)^T\bar{u}(k) \right),
\end{equation*}
and
$$
\bar{K}_\rho^*(P)=\arg\min_{K\in\K} \bar{J}_P (K,\rho).
$$
Using Lemma~\ref{prop:0}, we know that $\bar{K}_\rho^*(P)$ uniquely exists. We can find $\bar{J}_P (\bar{K}_\rho^*(P),\rho)$ using $X(\rho)$ as the unique positive definite solution of the discrete algebraic Riccati equation
\begin{equation} \label{eqn_Riccati}
\begin{split}
\tilde{A}^TX(\rho)\tilde{B}(\rho I+\tilde{B}^T&X(\rho)\tilde{B})^{-1}\tilde{B}^TX(\rho)\tilde{A}-\tilde{A}^TX(\rho)\tilde{A}+X(\rho)-I=0,
\end{split}
\end{equation}
with $\tilde{A}$ and $\tilde{B}$ defined in~(\ref{ABtilde}). According to~\cite{Komaroff1994}, we have
\begin{equation*}
\begin{split}
X(\rho) &\geq \tilde{A}^T(X_1^{-1}+(1/\rho)\tilde{B}\tilde{B}^T )^{-1}\tilde{A}+I \\ &=\tilde{A}^T(X_1-X_1\tilde{B}(\rho I+\tilde{B}^TX_1\tilde{B})^{-1}\tilde{B}^TX_1)\tilde{A}+I,
\end{split}
\end{equation*}
where
\begin{equation*}
\begin{split}
X_1&=\tilde{A}^T(I+(1/\rho)\tilde{B}\tilde{B}^T)^{-1}\tilde{A}+I.
\end{split}
\end{equation*}
Basic algebraic calculations show that
\begin{equation*}
\begin{split}
\lim_{\rho\rightarrow 0^+} X_1-X_1\tilde{B}&(\rho I+\tilde{B}^TX_1\tilde{B})^{-1}\tilde{B}^TX_1=\matrix{cc}{W & 0 \\ 0 & 0 }
\end{split}
\end{equation*}
where $W$ is defined in~(\ref{eqn:W}). According to~\cite{Kondo1986}, we know
$$
\lim_{\rho\rightarrow 0^+}\bar{J}_P (\bar{K}_\rho^*(P),\rho)=J_P(K^*_C(P))
$$
and as a result
\begin{equation*}
\begin{split}
X\hspace{-0.04in}=\hspace{-0.04in}\lim_{\rho\rightarrow 0^+}X(\rho) \geq \matrix{cc}{A & B \\0 & D}^T\hspace{-0.04in}\matrix{cc}{W & 0 \\ 0 & 0 }\hspace{-0.04in}\matrix{cc}{A & B \\0 & D}+I.
\end{split}
\end{equation*}
Equivalently, we get
\begin{equation} \label{eqn_16}
\matrix{cc}{X_{11} & X_{12} \\ X_{12}^T & X_{22} }\geq \matrix{cc}{A^TWA+I & A^TWB \\ BWA & BWB+I}.
\end{equation}
Now, we can calculate the cost of the optimal control design strategy as
\begin{equation} \label{eqn:cost_opt}
J_P(K^*_C(P))=\matrix{c}{x_0 \\ \xi(0)}^T \matrix{cc}{X_{11} & X_{12} \\ X_{12}^T & X_{22} } \matrix{c}{x_0 \\ \xi(0)}
\end{equation}
where
\begin{equation} \label{xi0}
\xi(0)=G_2 B^{-1}x_0+w_0=-(X_{22}^{-1}X_{12}^T+DB^{-1})x_0+w_0.
\end{equation}
If we put~(\ref{xi0}) in~(\ref{eqn:cost_opt}) and use the sub-Riccati equation
$$
X_{22}-I=BX_{11}B-BX_{12}X_{22}^{-1}X_{12}^TB,
$$
that is extracted from the Riccati equation in~(\ref{eqn_Riccati}) when $\rho=0$, we can simplify $J_P(K^*(P))$ in~(\ref{eqn:cost_opt}) to
\begin{small}
\begin{equation} \label{eqn_inter}
\begin{split}
&\matrix{c}{x_0 \\ -(X_{22}^{-1}X_{12}^T+DB^{-1})x_0+w_0}^T \hspace{-.07in}\matrix{cc}{X_{11} & X_{12} \\ X_{12}^T & X_{22} } \matrix{c}{x_0 \\ -(X_{22}^{-1}X_{12}^T+DB^{-1})x_0+w_0}\\&\hspace{.2in}=\matrix{c}{x_0 \\ w_0}^T \matrix{cc}{X_{11}-X_{12}X_{22}^{-1}X_{12}^T+B^{-1}DX_{22}DB^{-1} & -B^{-1}DX_{22} \\ -X_{22}DB^{-1} & X_{22} } \matrix{c}{x_0 \\ w_0}\\&\hspace{.2in}=\matrix{c}{x_0 \\ w_0}^T \matrix{cc}{B^{-1}(X_{22}+DX_{22}D-I)B^{-1} & -B^{-1}DX_{22} \\ -X_{22}DB^{-1} & X_{22} } \matrix{c}{x_0 \\ w_0}.
\end{split}
\end{equation}
\end{small}
Now, using~(\ref{eqn_16}) it is evident that $X_{22}\geq BWB+I$, and as a result
\begin{equation*}
\begin{split}
J_P(K^*_C(P))&\geq \matrix{c}{x_0 \\ w_0}^T \matrix{cc}{V_{11} & V_{12}B \\ BV_{12}^T & BV_{22}B } \matrix{c}{x_0 \\ w_0}.
\end{split}
\end{equation*}
where $V_{11}$, $V_{12}$, and $V_{22}$ are introduced in~(\ref{eqn:V1})-(\ref{eqn:V3}). The rest is only a straight forward matrix manipulation (factoring the matrix $B$).
\end{IEEEproof}

\subsection{Deadbeat Control Design Strategy}
In this subsection, we introduce the deadbeat control design strategy and give a useful lemma about its competitive ratio.

\begin{definition} \label{def:1} The deadbeat control design strategy $\Gamma^\Delta: \A(S_\p) \times \B(\epsilon_b) \times \D  \rightarrow \K$ is defined as
$$
\Gamma^\Delta(A,B,D)=\left[\begin{array}{c|c} D & -B^{-1}D^2 \\ \hline I & -B^{-1}(A+D) \end{array}\right].
$$
Using this control design strategy, irrespective of the value of the initial state $x_0$ and the initial disturbance $w_0$, the closed-loop system reaches the origin just in two time-steps. Note that the deadbeat control design strategy is a limited model information control design method since
$$
\Gamma_i^\Delta(A,B,D)=-(z-d_{ii})^{-1}b_{ii}^{-1}d_{ii}^2e_i^T-b_{ii}^{-1}(A_i+D_i)
$$
for each $1\leq i\leq n$. The cost of the deadbeat control design strategy $\Gamma^\Delta$ for any $P=(A,B,D,x_0,w_0)\in \p$ is
\begin{equation*}
J_P (\Gamma^\Delta(A,B,D))=\matrix{c}{x_0 \\ Bw_0}^T \matrix{cc}{Q_{11} & Q_{12}  \\ Q_{12}^T & Q_{22}} \matrix{c}{x_0 \\ Bw_0},
\end{equation*}
where
\begin{eqnarray}
&&\hspace{-.5in}Q_{11}=I+D^2(I+B^{-2})+A^TB^{-2}A \label{eqn_15_1} +DA^TB^{-2}AD+A^TB^{-2}D+DB^{-2}A, \\ &&\hspace{-.5in}Q_{12}=-D-A^TB^{-2}-DB^{-2}-DA^TB^{-2}A, \label{eqn_15_2} \\ &&\hspace{-.5in}Q_{22}=A^TB^{-2}A+B^{-2}+I. \label{eqn_15_3}
\end{eqnarray}
The closed-loop system with deadbeat control design strategy is shown in Figure~\ref{figure1}(\textit{a}). This feedback loop can be re-arranged as the one in Figure~\ref{figure1}(\textit{b}) which has two separate components. One component is a static-deadbeat control design strategy~\cite{Farokhi_ACC_2011} for regulating the state of the plant and the other one is the deadbeat observer for canceling the disturbance effect.
\end{definition}

\begin{lemma} \label{tho:1} Let the plant graph $G_{\p}$ contain no isolated node and $G_{\K}\supseteq G_\p$. Then, the competitive ratio of the deadbeat control design method $\Gamma^\Delta$ satisfies $r_\p(\Gamma^\Delta) \leq (2\epsilon_b^2+1+\sqrt{4\epsilon_b^2+1})/(2\epsilon_b^2)$.
\end{lemma}

\begin{IEEEproof} First, let us define the set of all real numbers that are greater than or equal to $r_\p(\Gamma^\Delta)$ as
$$
\mathcal{M}=\left\{\bar{\beta}\in\mathbb{R}\;\left|\; \frac{J_P(\Gamma^\Delta(A,B,D))}{J_P(K^*(P))}\leq \bar{\beta} \right.\;\forall P\in\p \right\}.
$$
It is evident that $J_P(K^*_C(P)) \leq J_P(K^*(P))$ for each plant $P\in\p$, irrespective of the control graph $G_\K$, and as a result
\begin{equation} \label{eqn:ineqKp*andKP*C}
\frac{J_P(\Gamma^\Delta(A,B,D))}{J_P(K^*(P))} \leq \frac{J_P(\Gamma^\Delta(A,B,D))}{J_P(K^*_C(P))}.
\end{equation}
Using Equation~(\ref{eqn:ineqKp*andKP*C}), Definition~\ref{def:1}, and Lemma~\ref{lem:1}, we get that $\beta$ belongs to the set $\mathcal{M}$ if
\begin{equation} \label{eqn_first_condition}
\begin{split}
\frac{\matrix{c}{x_0 \\ Bw_0 }^T \matrix{cc}{Q_{11} & Q_{12}  \\ Q_{12}^T & Q_{22}} \matrix{c}{x_0 \\ Bw_0 }}{\matrix{c}{x_0 \\ Bw_0 }^T \matrix{cc}{V_{11} & V_{12}  \\ V_{12}^T & V_{22}}  \matrix{c}{x_0 \\ Bw_0 }}\leq \beta,
\end{split}
\end{equation}
for all $A\in \A(S_\p)$, $B\in \B(\epsilon_b)$, $D\in \D $, $x_0\in \mathbb{R}^n$, and $w_0\in \mathbb{R}^n$ where $Q_{11}$, $Q_{12}$, and $Q_{22}$ are defined in~(\ref{eqn_15_1})-(\ref{eqn_15_3}) and $V_{11}$, $V_{12}$, and $V_{22}$ are defined in~(\ref{eqn:V1})-(\ref{eqn:V3}). The condition in~(\ref{eqn_first_condition}) is satisfied if and only if
\begin{equation*}
\begin{split}
\matrix{cc}{\beta V_{11}-Q_{11} & \beta V_{12}-Q_{12} \\ \beta V_{12}^T-Q_{12}^T & \beta V_{22}-Q_{22} } \geq 0,
\end{split}
\end{equation*}
for all $A\in \A(S_\p)$, $B\in \B(\epsilon_b)$, and $D\in \D $. Now, using Schur complement~\cite{Zhang2005}, we can show that $\beta$ belongs to the set $\mathcal{M}$ if both conditions
\begin{equation} \label{eqn:Z}
\begin{split}
Z&=\beta V_{22}-Q_{22}\\&=A^T(\beta (I+B^2)^{-1}-B^{-2})A+(\beta-1)(B^{-2}+I)\geq 0,
\end{split}
\end{equation}
and
\begin{equation} \label{eqn:condition1}
\begin{split}
\beta V_{11}-Q_{11}-&[\beta V_{12}-Q_{12}] [\beta V_{22}-Q_{22}]^{-1}[\beta V_{12}^T-Q_{12}^T]\geq 0,
\end{split}
\end{equation}
are satisfied for all matrices $A\in\A(S_\p)$, $B\in\B(\epsilon_b)$, and $D\in\D$. We can go further and simplify the condition in~(\ref{eqn:condition1}) to
\begin{equation} \label{eqn:condition1.5}
\begin{split}
\beta(W+DWD&+D^2B^{-2})-Q_{11}\\ &-\left[-DZ+A^TB^{-2} \right] Z^{-1}\left[-ZD+B^{-2}A\right]\geq 0,
\end{split}
\end{equation}
where $Z$ is introduced in~(\ref{eqn:Z}). For all $\beta \geq 1+1/\epsilon_b^2$, we know that $Z\geq (\beta-1)(B^{-2}+I)\geq 0$ and, as a result the condition
\begin{equation} \label{eqn:condition3}
\begin{split}
(\beta-1) I + &A^T\left(\beta(I+B^2)^{-1}-B^{-2}\right. \\& \hspace{-.2in} \left. -(\beta-1)^{-1} B^{-2}(B^{-2}+I)^{-1}B^{-2}\right)A \geq 0
\end{split}
\end{equation}
becomes a sufficient condition for the condition in~(\ref{eqn:condition1.5}) to be satisfied. Consequently, $\beta$ belongs to the set $\mathcal{M}$, if it is greater than or equal to $1+1/\epsilon_b^2$ and it satisfies the condition in~(\ref{eqn:condition3}). Thus, we get
$$
\left\{\beta\;|\;\beta \geq (2\epsilon_b^2+1+\sqrt{4\epsilon_b^2+1})/(2\epsilon_b^2) \right\}\subseteq \mathcal{M}.
$$
This concludes the proof.
\end{IEEEproof}

\jpfig{0.6}{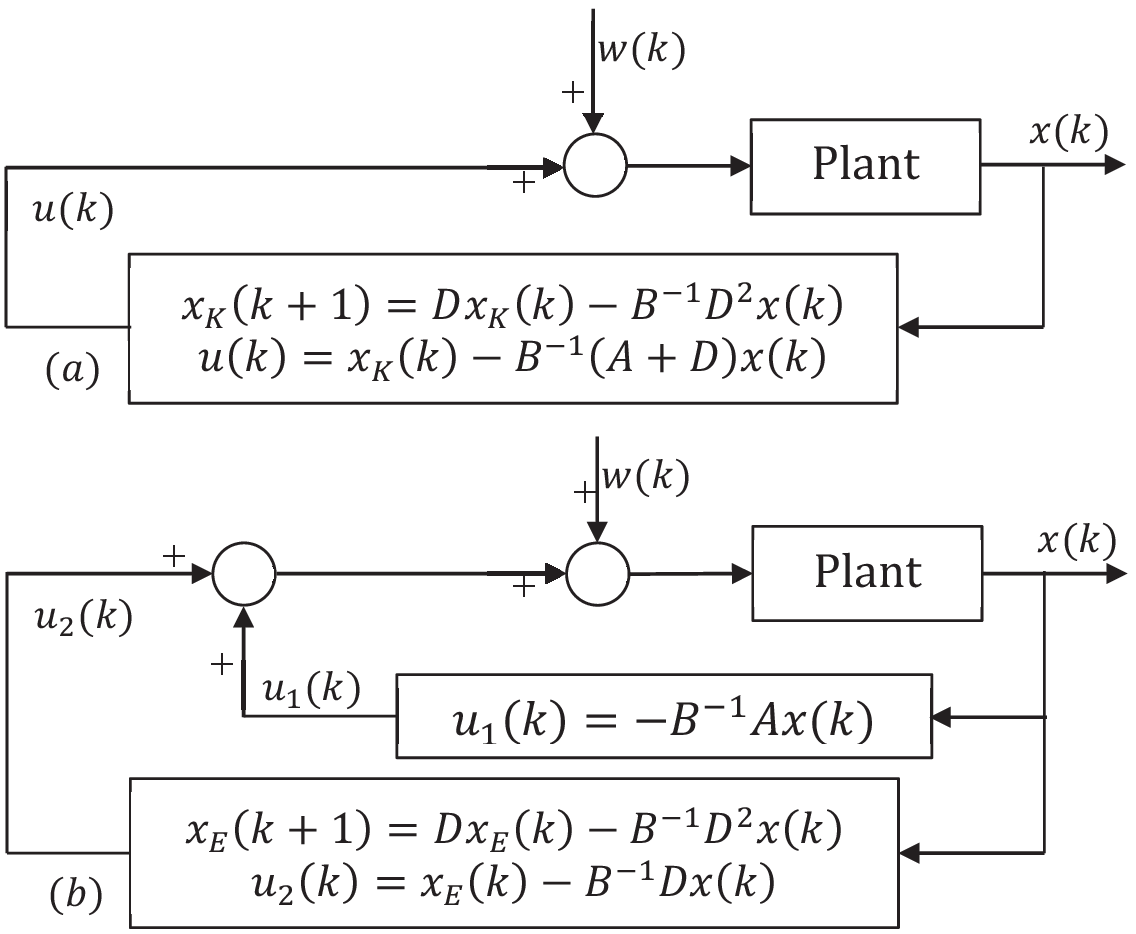}{figure1}{The closed-loop system with (\textit{a})~the deadbeat control design strategy $\Gamma^\Delta$ and (\textit{b})~rearranging this control design strategy as a static deadbeat control design and a deadbeat observer design.}

\section{Plant Graph Influence on Achievable Performance} \label{sec_3}
First, we need to give the following lemmas to make proof easier.

\begin{lemma} \label{prop:2} Let the plant graph $G_{\p}$ contain no isolated node and $G_{\K}\supseteq G_\p$. Let $P=(A,B,D,x_0,w_0)\in \p$ be a plant such that $A$ is a nilpotent matrix of degree two. Then, $J_P(K^*(P))=J_P(K^*_C(P))$.
\end{lemma}

\begin{IEEEproof} When matrix $A$ is nilpotent, based on the unique positive-definite solution of the discrete algebraic Riccati equation in~(\ref{eqn_Riccati}) when $\rho=0$, the optimal centralized controller $K^*_C(P)$ becomes
$$
K^*_C(P)= \left[\begin{array}{c|c} D & D(I+B^2)^{-1}B^{-1}A-B^{-1}D^2 \\ \hline I  & -(I+B^2)^{-1}BA-B^{-1}D \end{array}\right].
$$
Thus, $K^*_C(P)\in \K(S_\K)$ because $G_\K\supseteq G_\p$. Now, because $K^*(P)$ is the global optimal decentralized controller, it has a lower cost than any other decentralized controller $K\in \K(S_\K)$, and in particular
\begin{equation} \label{eqn:ineq1}
J_P(K^*(P))\leq J_P(K^*_C(P)).
\end{equation}
On the other hand, it is evident that
\begin{equation} \label{eqn:ineq2}
J_P(K^*_C(P)) \leq J_P(K^*(P)).
\end{equation}
The rest of the proof is a direct use of~(\ref{eqn:ineq1}) and~(\ref{eqn:ineq2}) simultaneously.
\end{IEEEproof}

\begin{lemma} \label{prop:1} Fix real numbers $a\in \mathbb{R}$ and $b\in \mathbb{R}$. For any $x\in \mathbb{R}$, we have $x^2+(a+bx)^2\geq a^2/(1+b^2)$. \end{lemma}

\begin{IEEEproof} Consider the function $x\mapsto x^2+(a+bx)^2$. Since this function is both continuously differentiable and strictly convex, we can find its unique minimizer as $\bar{x}=-ab/(1+b^2)$ by putting its derivative equal to zero. As a result, we get $x^2+(a+bx)^2\geq \bar{x}^2+(a+b\bar{x})^2= a^2/(1+b^2)$.
\end{IEEEproof}

\begin{lemma} \label{lem:2} Let the plant graph $G_{\p}$ contain no isolated node, the design graph $G_{\mathcal{C}}$ be a totally disconnected graph, and $G_\K\supseteq G_\p$. Furthermore, assume that node $i$ is not a sink in the plant graph $G_{\p}$. Then, the competitive ratio of control design strategy $\Gamma\in \comm$ is bounded only if $a_{ij}+b_{ii}(d_\Gamma)_{ij}(A,B,D)=0$ for all $j\neq i$ and all matrices $A\in \A(S_\p)$, $B\in \B(\epsilon_b)$, and $D\in \D$. \end{lemma}

\begin{IEEEproof} The proof is by contrapositive. Assume that the matrices $\bar{A}\in \A(S_\p)$, $B\in \B(\epsilon_b)$, $D\in \D $, and indices $i$ and $j$ exist such that $i\neq j$ and $\bar{a}_{ij}+b_{ii} (d_\Gamma)_{ij}(\bar{A},B,D)\neq 0$ for some control design strategy $\Gamma\in\comm$. Let $1\leq \ell\leq n$ be an index such that $\ell \neq i$ and $(s_\p)_{\ell i}\neq 0$ (such an index exists because node $i$ is not a sink in the plant graph). Define matrix $A$ such that $A_i=\bar{A}_i$, $A_\ell=re_i^T$, and $A_t=0$ for all $t\neq i,\ell$. It is evident that $\Gamma_i(\bar{A},B,D)=\Gamma_i(A,B,D)$ since the design graph is a totally disconnected graph. Using the structure of the cost function in~(\ref{cost_easy}) and plant dynamics in~(\ref{eqn_1}), the cost of the control design strategy $\Gamma$ when $w_0=e_j$ and $x_0=0$ satisfies
\begin{equation*}
\begin{split}
J_{P}(\Gamma(A,B,D))&\geq (u_\ell(2)+w_\ell(2))^2+x_\ell(3)^2 \\& =(u_\ell(2)+w_\ell(2))^2+(rx_i(2)+b_{\ell\ell}(u_\ell(2)+w_\ell(2)))^2.
\end{split}
\end{equation*}
With the help of Lemma~\ref{prop:1} and the fact that $x_i(2)=(a_{ij}+b_{ii}(d_\Gamma)_{ij}(A,B,D))b_{jj}$ (see~Figure~\ref{figure2}), we get
\begin{equation*}
\begin{split}
J_{P}(\Gamma(A,B,D))&\geq r^2x_i(2)^2/(1+b_{\ell\ell}^2) \\&= (a_{ij}+b_{ii}(d_\Gamma)_{ij}(A,B,D))^2b_{jj}^2 r^2/(1+b_{\ell\ell}^2).
\end{split}
\end{equation*}
The cost of the deadbeat control design strategy is
\begin{equation*}
\begin{split}
J_{P}(\Gamma^\Delta(A,B,D))&=e_j^TB^T(A^TB^{-2}A+B^{-2}+I)Be_j\\&=b_{jj}^2+1+a_{ij}^2b_{jj}^2/b_{ii}^2.
\end{split}
\end{equation*}
Using the inequality
\begin{equation*}
\begin{split}
r_\p(\Gamma)&=\sup_{P \in \p } \frac{J_P (\Gamma(A,B,D))}{J_P (K^*(P))}\\&= \sup_{P \in \p } \left[ \frac{J_P (\Gamma(A,B,D))}{J_P (\Gamma^\Delta(A,B,D))} \frac{J_P (\Gamma^\Delta(A,B,D))}{J_P (K^*(P))} \right] \\& \geq \sup_{P \in \p } \frac{J_P (\Gamma(A,B,D))}{J_P (\Gamma^\Delta(A,B,D))},
\end{split}
\end{equation*}
gives
\begin{equation*}
\begin{split}
r_\p(\Gamma)\geq \frac{(a_{ij}+b_{ii}(d_\Gamma)_{ij}(A,B,D))^2b_{jj}^2}{(1+b_{\ell\ell}^2)(b_{jj}^2+1+a_{ij}^2b_{jj}^2/b_{ii}^2)}\lim_{r\rightarrow \infty}r^2 =\infty.
\end{split}
\end{equation*}
This proves the statement by contrapositive.
\end{IEEEproof}

Now, we are ready to tackle the problem~(\ref{eqn:0}). As the main results of the paper crucially depends on the properties of the plant graph, we split these results to two different subsections.

\subsection{Plant Graphs without Sinks}
In this section, we assume that there is no sink in the plant graph, and we try to find the best control design strategy in terms of the competitive ratio and the domination.

\begin{theorem} \label{tho:2}  Let the plant graph $G_{\p}$ contain no isolated node and no sink, the design graph $G_{\mathcal{C}}$ be a totally disconnected graph, and $G_\K\supseteq G_\p$. Then, the following statements hold:
\newline (\textit{a}) The competitive ratio of any control design strategy $\Gamma\in\comm$ satisfies $r_{\p} (\Gamma)\geq r_{\p} (\Gamma^\Delta) = (2\epsilon_b^2+1+\sqrt{4\epsilon_b^2+1})/(2\epsilon_b^2)$. \newline (\textit{b}) The control design strategy $\Gamma^\Delta$ is undominated, if and only if, there is no sink in the plant graph $G_{\p}$.
\end{theorem}

\begin{IEEEproof} First, let us prove statement~(\textit{a}). It is always possible to pick indices $j\neq i$ such that $(s_\p)_{ji}\neq 0$ since there is no isolated node in the plant graph $G_\p$. Let us define a one-parameter family of matrices $\{A(r)\}$ where $A(r)=re_je_i^T$ for each $r\in \mathbb{R}$. In addition, let $B=\epsilon_b I$ and $D=I$. According to Lemma~\ref{lem:2}, $r_\p(\Gamma)$ is bounded only if $r+\epsilon_b(d_\Gamma)_{ji}(r)=0$. Therefore, there is no loss of generality in assuming that $(d_\Gamma)_{ji}(r)=-r/\epsilon_b$ because otherwise $r_\p(\Gamma)$ is infinity and the inequality $r_\p(\Gamma)\geq r_\p(\Gamma^\Delta)$ is trivially satisfied (considering that using Lemma~\ref{tho:1} we know $r_\p(\Gamma^\Delta)$ is bounded). For each $r\in \mathbb{R}$, the matrix $A(r)$ is a nilpotent matrix of degree two. Thus, using Lemma~\ref{prop:2}, we get $J_P(K^*(P))=J_P(K^*_C(P))$ for this special plant. The unique positive definite solution of the discrete algebraic Riccati equation in~(\ref{eqn_Riccati}) for a fixed $r$ (when $\rho=0$) is
$$
X=\matrix{cc}{A(r)^TA(r) & \epsilon_b A(r)^T \\ \epsilon_b A(r) & \epsilon_b^2/(1+\epsilon_b^2) A(r)^TA(r) + \epsilon_b^2I}+I.
$$
Thus, the cost of the optimal control design strategy for
\begin{equation} \label{eqn_x_0}
x_0=\frac{(\epsilon_b^2 + 1)(\sqrt{4\epsilon_b^2 + 1}+ 1)}{2\epsilon_br} e_i,
\end{equation}
and
\begin{equation} \label{eqn_w_0}
w_0=\frac{(\epsilon_b^2 + 1)(\sqrt{4\epsilon_b^2 + 1} + 1)}{2\epsilon_b^2r}e_i-e_j,
\end{equation}
is equal to
\begin{equation*}
\begin{split}
J_P(K^*(P))=& \frac{\epsilon_b^2\sqrt{4\epsilon_b^2 + 1} + 5\epsilon_b^2 + 4\epsilon_b^4 + \sqrt{4\epsilon_b^2 + 1} + 1}{2\epsilon_b^2} \\ &\hspace{0.2in} + \frac{(2\epsilon_b^2 + \sqrt{4\epsilon_b^2 + 1} + 1)\sqrt{4\epsilon_b^2 + 1}}{2\epsilon_b^2r^2},
\end{split}
\end{equation*}
On the other hand, for each $r\in \mathbb{R}$, the cost of the control design strategy $\Gamma$ for $x_0$ and $w_0$ given in~(\ref{eqn_x_0})~and~(\ref{eqn_w_0}) is lower-bounded by
\begin{equation*}
\begin{split}
J_P(\Gamma(A,B,D)) &\geq (u_j(0)+w_j(0))^2 + x_j(1)^2
\\& =\frac{(\epsilon_b^2 + 1)(3\epsilon_b^2\sqrt{4\epsilon_b^2 + 1} + 5\epsilon_b^2 + 4\epsilon_b^4 + \sqrt{4\epsilon_b^2 + 1} + 1)}{2\epsilon_b^4}.
\end{split}
\end{equation*}
Therefore, for any $\Gamma\in\comm$, we have
\begin{equation} \label{eqn:ineqGamma}
\begin{split}
r_\p(\Gamma)\geq \lim_{r\rightarrow \infty} \frac{J_P(\Gamma(A,B,D))}{J_P(K^*(P))}=\frac{2\epsilon_b^2+1+\sqrt{4\epsilon_b^2+1}}{2\epsilon_b^2}.
\end{split}
\end{equation}
Considering the fact that $\Gamma^\Delta$ also belongs to $\comm$, the rest is a simple combination of~(\ref{eqn:ineqGamma}) and Lemma~\ref{tho:1}.
\par Now, we can prove statement~(\textit{b}). The ``if'' part of the proof is done by constructing plants $P=(A,B,D,x_0,w_0)\in \p$ that satisfy $J_P(\Gamma(A,B,D))>J_P(\Gamma^\Delta(A,B,D))$ for any control design method $\Gamma\in \comm \setminus \{\Gamma^\Delta \}$. For the ``only if'' part, we show that $\Gamma^{\Theta}$ introduced later in~(\ref{eqn:11}) dominates $\Gamma^\Delta$ when $G_{\mathcal{P}}$ has at least one sink. See~\cite[p.124]{Farokhi_Technical_2011} for the detailed proof.
\end{IEEEproof}

\jpfig{0.8}{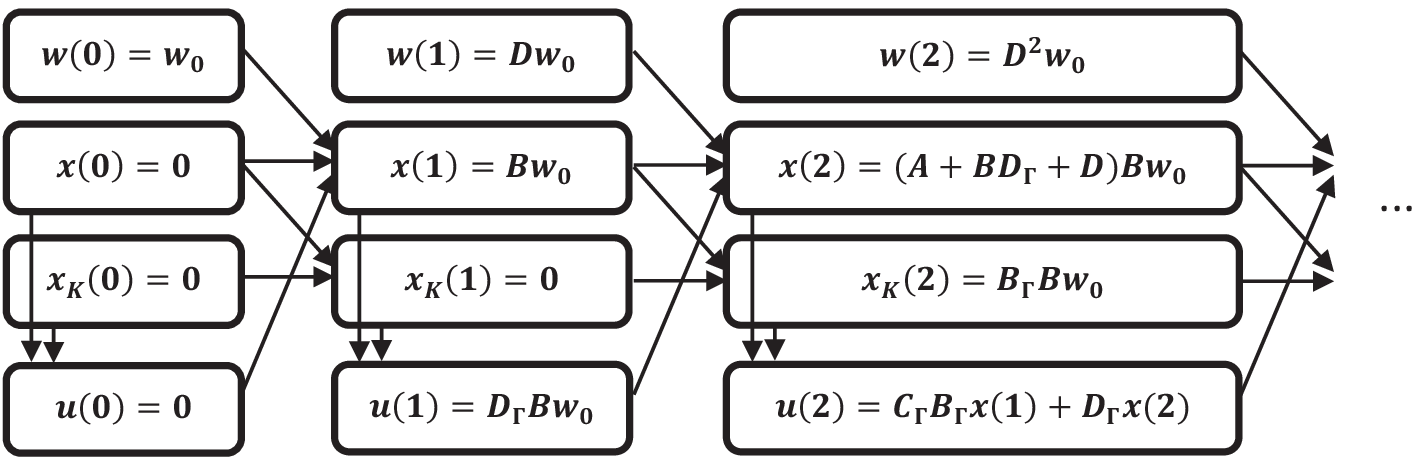}{figure2}{State evolution of the closed-loop system when $x_0=0$.}

Theorem~\ref{tho:2} shows that the deadbeat control design method $\Gamma^\Delta$ is an undominated minimizer of the competitive ratio $r_{\p}$ over the set of limited model information design methods $\comm$.

\subsection{Plant Graphs with Sinks}
In this section, we study the case where there are $c\geq 1$ sinks in the plant graph $G_\p$. By renumbering the sinks as subsystems number $n-c+1,\dots,n$, the matrix $S_\p$ can be written as
\begin{equation} \label{eqn:10}
S_{\p} = \matrix{c|c}{ (S_{\p})_{11} & 0_{(q-c)\times (c)} \\ \hline (S_{\p})_{21} & (S_{\p})_{22} },
\end{equation}
where
\begin{equation*}
(S_\p)_{11}=\matrix{ccc}{ (s_{\p})_{11} & \cdots & (s_{\p})_{1,n-c}\\ \vdots & \ddots & \vdots \\ (s_{\p})_{n-c,1} & \cdots & (s_{\p})_{n-c,n-c} },
\end{equation*}
\begin{equation*}
(S_\p)_{21}=\matrix{ccc}{ (s_{\p})_{n-c+1,1} & \cdots & (s_{\p})_{n-c+1,n-c}\\ \vdots & \ddots & \vdots \\ (s_{\p})_{n,1} & \cdots & (s_{\p})_{n,n-c} },
\end{equation*}
and $(S_\p)_{22}=\diag \left((s_{\p})_{n-c+1,n-c+1},\dots,(s_{\p})_{nn}\right)$. From now on, without loss of generality, we assume that the structure matrix is the one defined in~(\ref{eqn:10}). For all plants $P\in\p$, control design method $\Gamma^\Theta$ is defined as
\begin{equation} \label{eqn:11}
\Gamma^\Theta(A,B,D)=\left[\begin{array}{c|c} D & B^{-1}DF(A,B)A-B^{-1}D^2 \\ \hline I  & B^{-1}(F(A,B)-I)A-B^{-1}D \end{array}\right]
\end{equation}
where
\begin{equation*}
F(A,B)=\textrm{diag}\left(0,\dots,0,f_{n-c+1}(A,B),\dots,f_n(A,B)\right),
\end{equation*}
and
\begin{equation*}
f(A,B)=\frac{2}{b_{ii}^2+a_{ii}^2+1+
\sqrt{(a_{ii}^2+b_{ii}^2)^2+2(b_{ii}^2-a_{ii}^2)+1}}
\end{equation*}
for all $n-c+1\leq i\leq n$.

The control design strategy $\Gamma^\Theta$ applies the deadbeat to every subsystem that is not a sink and, for every sink, applies the same optimal control law as if the node were decoupled from the rest of the graph.

\begin{lemma} \label{tho:5}
Let the plant graph $G_{\p}$ contain no isolated node and at least one sink and $G_\K\supseteq G_\p$. Then, the competitive ratio of the design method $\Gamma^\Theta$ introduced in~(\ref{eqn:11}) is
$$
r_{\p}(\Gamma^\Theta)=\left\{ \begin{array}{ll} (2\epsilon_b^2+1+\sqrt{4\epsilon_b^2+1})/(2\epsilon_b^2), &  \mbox{if} \; (S_{\p})_{11} \mbox{ is not diagonal}, \\ 1, & \mbox{if} \; (S_{\p})_{11}=0 \;\&\; (S_{\p})_{22}=0. \end{array} \right.
$$
\end{lemma}

\begin{IEEEproof} Based the proof of the ``only if'' part of statement~(\textit{b}) of Theorem~\ref{tho:2}, we know that
\begin{equation*}
J_P(\Gamma^\Theta(A,B,D))\leq J_P(\Gamma^\Delta(A,B,D)),
\end{equation*}
for all $P=(A,B,D,x_0,w_0)\in\p$ and as a result
\begin{equation*}
\begin{split}
r_\p(\Gamma^\Theta)&=\sup_{P\in \p} \frac{J_{P}(\Gamma^\Theta(A,B,D))}{J_{P}(K^*(P))} \\&\leq \sup_{P\in \p} \frac{J_{P}(\Gamma^\Delta(A,B,D))}{J_{P}(K^*(P))}  \leq \frac{2\epsilon_b^2+1+\sqrt{4\epsilon_b^2+1}}{2\epsilon_b^2}.
\end{split}
\end{equation*}
Now if $(S_\p)_{11}$ has an off-diagonal entry, then there exist $1\leq i,j\leq n-c$ and $i\neq j$ such that $(s_\p)_{ji}\neq 0$. Using the second part of the proof of Theorem~\ref{tho:2}, it is easy to see
$$
r_\p(\Gamma^\Theta) \geq \frac{2\epsilon_b^2+1+\sqrt{4\epsilon_b^2+1}}{2\epsilon_b^2},
$$
because the control design $\Gamma^\Theta$ acts like the deadbeat control design strategy on that part of the system. Using both these inequalities proves the statement.
\par If $(S_\p)_{11}=0$ and $(S_\p)_{22}=0$, every matrix $A$ with structure matrix $S_\p$ becomes a nilpotent matrix of degree two. Thus, according to Lemma~\ref{prop:2}, we get that $J_P(K^*(P))=J_P(K^*_C(P))$, and based on the unique solution of the associated discrete algebraic Riccati equation, for this plant, the optimal centralized control design is
$$
K^*_C(P)=\left[\begin{array}{c|c} D & D(I+B^2)^{-1}B^{-1}A-B^{-1}D^2 \\ \hline I  & -(I+B^2)^{-1}BA-B^{-1}D \end{array}\right],
$$
which is exactly equal to $\Gamma^\Theta(A,B,D)$. Thus, $r_\p(\Gamma^\Theta)=1$.
\end{IEEEproof}

\begin{theorem} \label{tho:6} Let the plant graph $G_{\p}$ contain no isolated node and contain at least one sink, the design graph $G_{\mathcal{C}}$ be a totally disconnected graph, and $G_{\K}\supseteq G_\p$. Then, the following statements hold:
\newline (\textit{a}) The competitive ratio of any control design strategy $\Gamma\in\comm$ satisfies $ r_{\p}(\Gamma)\geq (2\epsilon_b^2+1+\sqrt{4\epsilon_b^2+1})/(2\epsilon_b^2)$, if $(S_{\p})_{11}$ is not diagonal.
\newline (\textit{b}) The control design method $\Gamma^\Theta$ is undominated by all limited model information control design methods in $\comm$.
\end{theorem}

\begin{IEEEproof} First, we prove statement~(\textit{a}). Suppose that $(S_\p)_{11}\neq 0$ and $(S_\p)_{11}$ is not a diagonal matrix, then there exist $1\leq i,j \leq n-c$ and $i\neq j$ such that $(s_\p)_{ji}\neq 0$. Consider the family of matrices $A(r)$ defined by $A(r)=re_je_i^T$. Based on Lemma~\ref{lem:2}, if we want to have a bounded competitive ratio, the control design strategy should satisfy $r+b_{jj}(d_\Gamma)_{ji}(A(r),B,D)=0$ (because node $1\leq j\leq n-c$ is not a sink). The rest of the proof is similar to the proof of Theorem~\ref{tho:2}.
\par See~\cite[p.130]{Farokhi_Technical_2011} for the detailed proof of statement~(\textit{b}).
\end{IEEEproof}

Combining Lemma~\ref{tho:5}~and~Theorem~\ref{tho:6} illustrates that if $(S_{\p})_{11}\neq 0$ is not diagonal, the control design method $\Gamma^\Theta$ has the smallest ratio achievable by limited model information control methods. Thus, it is a solution to the problem~(\ref{eqn:0}). Furthermore, if $(S_{\p})_{11}$ and $(S_{\p})_{22}$ are both zero, then $\Gamma^\Theta$ becomes equal to $K^*$. This shows that $\Gamma^\Theta$ is a solution to the problem~(\ref{eqn:0}) in this case too. The rest of the cases are still open.

\section{Design Graph Influence on Achievable Performance}  \label{sec:Gc}
In the previous section, we solved the optimal control design under limited model information when $G_{\mathcal{C}}$ is a totally disconnected graph. In this section, we study the necessary amount of information needed in each subsystem to ensure the existence of a limited model information control design strategy with a better competitive ratio than $\Gamma^{\Delta}$ and $\Gamma^{\Theta}$.

\begin{theorem} \label{tho:8} Let the plant graph $G_\p$ and the design graph $G_\comm$ be given and $G_{\K}\supseteq G_\p$. Then, we have $r_{\p}(\Gamma) \geq (2\epsilon_b^2+1 +\sqrt{4\epsilon_b^2+1})/(2\epsilon_b^2)$ for all $\Gamma \in \mathcal{C}$ if $G_{\mathcal{P}}$ contains the path $i \rightarrow j \rightarrow \ell$ with distinct nodes $i$, $j$, and $\ell$ while $(\ell,j)\notin E_\comm$.
\end{theorem}

\begin{IEEEproof} See~\cite[p.132]{Farokhi_Technical_2011} for the detailed proof. \end{IEEEproof}

\section{Conclusions} \label{sec_4}
We studied the design of optimal dynamic disturbance accommodation controllers under limited plant model information. To do so, we investigated the relationship between closed-loop performance and the control design strategies with limited model information using the performance metric called the competitive ratio. We found an explicit minimizer of the competitive ratio and showed that this minimizer is also undominated. Possible future work will focus on extending the present framework to situations where the subsystems are not scalar.

\bibliography{ref}

\begin{thebibliography}{10}

\bibitem{Giulietti2000}
F.~Giulietti, L.~Pollini, and M.~Innocenti, ``Autonomous formation flight,''
  {\em Control Systems Magazine, IEEE}, vol.~20, no.~6, pp.~34 -- 44, 2000.

\bibitem{Swaroop1999}
D.~Swaroop and J.~K. Hedrick, ``Constant spacing strategies for platooning in
  automated highway systems,'' {\em Journal of Dynamic Systems, Measurement,
  and Control}, vol.~121, no.~3, pp.~462--470, 1999.

\bibitem{Dunbar2007}
W.~Dunbar, ``Distributed receding horizon control of dynamically coupled
  nonlinear systems,'' {\em Automatic Control, IEEE Transactions on}, vol.~52,
  no.~7, pp.~1249 --1263, 2007.

\bibitem{Negenborn2010}
R.~R. Negenborn, Z.~Lukszo, and H.~Hellendoorn, eds., {\em Intelligent
  Infrastructures}, vol.~42.
\newblock Springer, 2010.

\bibitem{Levine1971}
W.~Levine, T.~Johnson, and M.~Athans, ``Optimal limited state variable feedback
  controllers for linear systems,'' {\em Automatic Control, IEEE Transactions
  on}, vol.~16, no.~6, pp.~785 -- 793, 1971.

\bibitem{Rotkowitz2006}
M.~Rotkowitz and S.~Lall, ``A characterization of convex problems in
  decentralized control,'' {\em Automatic Control, IEEE Transactions on},
  vol.~51, no.~2, pp.~274 -- 286, 2006.

\bibitem{Voulgaris200351}
P.~G. Voulgaris, ``Optimal control of systems with delayed observation sharing
  patterns via input-output methods,'' {\em Systems \& Control Letters},
  vol.~50, no.~1, pp.~51 -- 64, 2003.

\bibitem{Farokhi_ACC_2011}
F.~Farokhi, C.~Langbort, and K.~H. Johansson, ``Control design with limited
  model information,'' in {\em American Control Conference, Proceedings of
  the}, pp.~4697 -- 4704, 2011.

\bibitem{Farokhi_Technical_2011}
F.~Farokhi, ``Decentralized control design with limited plant model
  information,'' Licentiate Thesis, 2012.
\newblock \url{http://urn.kb.se/resolve?urn=urn:nbn:se:kth:diva-63858}.

\bibitem{Johnson1968}
C.~Johnson, ``Optimal control of the linear regulator with constant
  disturbances,'' {\em Automatic Control, IEEE Transactions on}, vol.~13,
  no.~4, pp.~416 -- 421, 1968.

\bibitem{Anderson1971}
B.~D.~O. Anderson and J.~B. Moore, {\em {Linear Optimal Control}}.
\newblock Prentice-Hall, 1971.

\bibitem{Langbort2010}
C.~Langbort and J.-C. Delvenne, ``Distributed design methods for linear
  quadratic control and their limitations,'' {\em Automatic Control, IEEE
  Transactions on}, vol.~55, no.~9, pp.~2085 --2093, 2010.

\bibitem{Molinari1975}
B.~Molinari, ``The stabilizing solution of the discrete algebraic
  \textsc{R}iccati equation,'' {\em Automatic Control, IEEE Transactions on},
  vol.~20, no.~3, pp.~396 -- 399, 1975.

\bibitem{Komaroff1994}
N.~Komaroff, ``Iterative matrix bounds and computational solutions to the
  discrete algebraic \textsc{R}iccati equation,'' {\em Automatic Control, IEEE
  Transactions on}, vol.~39, no.~8, pp.~1676 --1678, 1994.

\bibitem{Kondo1986}
R.~Kondo and K.~Furuta, ``On the bilinear transformation of \textsc{R}iccati
  equations,'' {\em Automatic Control, IEEE Transactions on}, vol.~31, pp.~50
  -- 54, Jan. 1986.

\bibitem{Zhang2005}
F.~Zhang, {\em {The Schur Complement and Its Applications}}.
\newblock Springer, 2005.

\end{thebibliography}
\bibliographystyle{ieeetr}

\end{document}